\documentclass[preprint,11pt,times]{elsarticle}

\usepackage{setspace}
\usepackage{geometry}
\usepackage{url}
\usepackage{graphicx,graphics}
\usepackage{subfigure}
\usepackage{epsfig}
\usepackage{psfrag}
\usepackage{hyperref}
\usepackage{amsmath,amssymb,amsfonts,amsthm}
\usepackage{mathrsfs}
\usepackage{color}
\usepackage{float}
\usepackage{natbib}
\usepackage{tabularx}
\usepackage{multirow}

\usepackage{pgfplots}
\usetikzlibrary{plotmarks}

\theoremstyle{remark}

\newcommand{\Eref}[1]{Equation (\ref{#1})}
\newcommand{\fref}[1]{Figure (\ref{#1})}

\newcommand{\rmd}{\mathrm{d}}
\newcommand{\bveps}{\boldsymbol{\varepsilon}}
\newcommand{\bvsig}{\boldsymbol{\sigma}}

\newcommand{\bb}{\mathbf{b}}
\newcommand{\bigb}{\mathbf{B}}

\newcommand{\dd}{\mathbf{D}}

\newcommand{\ff}{\mathbf{F}}
\newcommand{\kk}{\mathbf{K}}
\newcommand{\bn}{\mathbf{N}}

\newcommand{\qq}{\mathbf{q}}

\newcommand{\ttt}{\mathbf{T}}
\newcommand{\uu}{\mathbf{u}}

\newcommand{\xx}{\mathbf{x}}

\begin{document}

\begin{frontmatter}

\title{A hybrid T-Trefftz polygonal finite element for linear elasticity}       

\author[india]{Kalyan Bhattacharjee}
\author[unsw]{Sundararajan Natarajan }
\author[lux]{St\'ephane Bordas \corref{cor1}\fnref{fn1} }

\cortext[cor1]{Corresponding author}

\address[india]{Department of Ocean Engineering \& Naval Architecture, Indian Institute of Technology, Kharagpur, India}
\address[unsw]{School of Civil \& Environmental Engineering, The University of New South Wales, Sydney, Australia}
\address[lux]{Faculté des Sciences, de la Technologie et de la Communication, University of Luxembourg, Luxembourg}

\fntext[fn1]{Faculté des Sciences, de la Technologie et de la Communication, University of Luxembourg, Luxembourg. Email: stephane.bordas@uni.lu}

\begin{abstract}

In this paper, we construct hybrid T-Trefftz polygonal finite elements. The displacement field within the polygon is represented by the homogeneous solution to the governing differential equation, also called as the T-complete set. On the boundary of the polygon, a conforming displacement field is independently defined to enforce continuity of the displacements across the element boundary. An optimal number of T-complete functions are chosen based on the number of nodes of the polygon and degrees of freedom per node. The stiffness matrix is computed by the hybrid formulation with auxiliary displacement frame. Results from the numerical studies presented for a few benchmark problems in the context of linear elasticity shows that the proposed method yield highly accurate results.
	
\end{abstract}


\end{frontmatter}


\section{Introduction}

\subsection{Background}
The finite element method (FEM) relies on discretizing the domain with non-overlapping regions, called the `elements'. In the conventional FEM, the topology of the elements were restricted to triangles and quadrilaterals in 2D or tetrahedrals and hexahedrals in 3D. The use of such standard shapes, simplifies the approach, however, this may require sophisticated (re-) meshing algorithms to either generate high-quality meshes or to capture the topological changes. Moreover, the accuracy of the solution depends on the quality of the element employed. Lee and Bathe~\cite{leebathe1993} observed that the shape functions lose their ability to reproduce the displacement fields when the mesh is distorted. In an effort to overcome the limitations of the FEM, the research has been focussed on:
\begin{itemize}
\item de-coupling geometry and analysis, for example, meshfree methods~\cite{gingoldmonaghan1977,belytschkolu1994}, PU enrichment~\cite{melenkbabuvska1996}, Immersed boundary method~\cite{peskin2002}.
\item Improve the element formulations
\begin{itemize}
\item Strain smoothing~\cite{liunguyen2007}
\item Unsymmetric formulation~\cite{rajendranooi2007,rajendran2010}
\item hybrid Trefftz FEM~\cite{szeliu2010,wangqin2011}
\item Polygonal FEM~\cite{sukumarmalsch2006}
\end{itemize}
\item coupling geometry and analysis, for example, isogeometric analysis~\cite{kaganfischer2003,hughescottrell2006}.
\end{itemize}

\begin{figure}[htpb]
\centering
\subfigure[]{\includegraphics[scale=0.5]{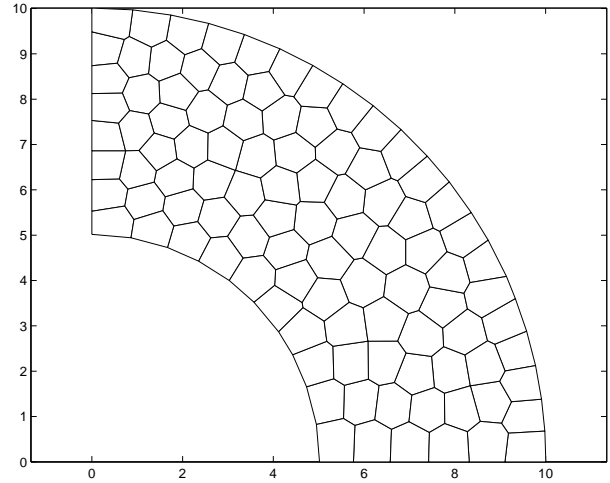}}
\subfigure[]{\includegraphics[scale=0.5]{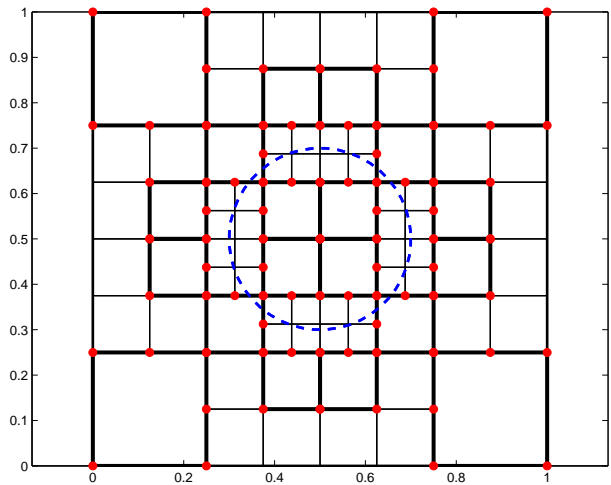}}
\caption{Domain discretized with finite elements: (a) discretization with arbitrary polygons and (b) adaptive refinement leading to a quadtree mesh, where the `dotted' line represents an inner boundary and the `highlighted' elements are the elements with hanging nodes.}
\label{fig:pmeshexample}
\end{figure}

In this study, we focus on the work on the improvement of finite element formulation, esp., the polygonal FEM. The use of elements with arbitrary number of sides provides flexibility in automatic mesh manipulation. For example, the domain can be discretized without a need to maintain a particular element topology (see \fref{fig:pmeshexample}). Moreover, this is advantageous in adaptive mesh refinement, where a straightforward subdivision of individual elements usually results in hanging nodes (see \fref{fig:pmeshexample}). Conventionally, this is eliminated by introducing additional edges/faces to retain conformity. This can be avoided if we can compute directly on polyhedral meshes with hanging nodes. 

In 1971~\cite{wachspress1971}, Wachspress developed a method based on rational basis functions for generalization to elements with arbitrary number of sides. However, these elements were not used because of the associated difficulties, such as in constructing the basis functions and the numerical integration of the basis functions. Thanks to advancement in mathematical softwares, viz., MATHEMATICA\textsuperscript{\textregistered} and MAPLE\textsuperscript{\textregistered} and the pioneering work of Alwood and Cornes~\cite{alwoodcornes1969}, Sukumar and Tabarraei~\cite{sukumartabarraei2004}, Dasgupta~\cite{dasgupta2003}, to name a few. Now discretization of the domain with finite elements having arbitrary number of sides has gained increasing attention~\cite{daveigabrezzi2013}. Once the basis functions are constructed, the conventional Galerkin procedure is normally employed to solve the governing equations over the polygonal/polyhedral meshes.

Ghosh \textit{et al.,}~\cite{moorthyghosh2000} developed the Vorono\"{i} cell finite element method (VCFEM) to model the mechanical response of heterogeneous microstructures of composites and porous materials with heterogeneities. The VCFEM is based on the assumed stress hybrid formulation and was further developed by Tiwary \textit{et al.,}~\cite{tiwaryhu2007} to study the behaviour of microstructures with irregular geometries. Rashid and Gullet~\cite{rashidgullet2000} proposed a variable element topology finite element method (VETFEM), in which the shape functions are constructed using a constrained minimization procedure. Sukumar~\cite{sukumar2003} used Vorono\"{i} cells and natural neighbor interpolants to develop a finite difference method on unstructured grids. Liu \textit{et al.,}~\cite{} generalized the concept of strain smoothing technique to arbitrarily shaped polygons. The main idea is to write the strain as the divergence of a spatial average of the compatible strain field. On another front, a fundamental solution less method was introduced by Wolf and Song~\cite{wolfsong2001}. It shares the advantages of the FEM and the boundary element method (BEM). Like the FEM, no fundamental solution is required and like the BEM, the spatial dimension is reduced by one, since only the boundary need to be discretized, resulting in a decrease in the total degrees of freedom. Ooi \textit{et al.,}~\cite{ooisong2012} employed scaled boundary formulation in polygonal elements to study crack propagation. Earlier, the SBFEM. Apart from the aforementioned formulations, recent studies, among others include developing polygonal elements based on the virtual nodes~\cite{tangwu2009} and the virtual element methods~\cite{veigabrezzi2013}. The other possible approach is to employ basis functions that satisfy the differential equation locally~\cite{qin2005,copelandlanger2009}. This method has been studied in detail in~\cite{hofreitherlanger2010,weisser2011} and extended to higher order polygons in~\cite{rjasanowwei$ss$er2012,weisser2012}. It is beyond the scope of this paper to review advances in polygonal finite element methods. There are different techniques to compute the basis functions over arbitrary polygons. Some of them include: (a) Using length and area measures~\cite{wachspress1971}; (b) Natural neighbour interpolants~\cite{sukumarmoran1998}; (c) Maximum entropy approximant~\cite{sukumar2013}; (d) Harmonic shape functions~\cite{bishop2013}. However, approximation functions on polygonal elements are usually non-polynomial, which introduces difficulties in the numerical integration. Improving numerical integration over polytopes has gained increasing attention~\cite{sukumartabarraei2004,natarajanbordas2009,mousavixiao2010,talischipaulino2013}. Apart from the aforementioned formulations, recent studies, among others include developing polygonal elements based on the virtual nodes~\cite{tangwu2009} and the virtual element methods~\cite{veigabrezzi2013}. It is beyond the scope of this paper to review advances in polygonal finite element methods. Interested readers are referred to the literature~\cite{friesmatthies2003,sukumarmalsch2006} and references therein for detailed discussion.

The other possible approach is to employ basis functions that satisfy the differential equation locally~\cite{qin2005,copelandlanger2009}. Zienkiewicz~\cite{zienkiewicz1997} presented a concise discussion on different approximation procedures to differential equations. It was shown that Trefftz type approximation is a particular form of weighted residual approximation. This can be used to generate \emph{hybrid finite elements}. Earlier studies employed boundary type approximation associated with Trefftz to develop special type finite elements, for example, elements with holes/voids~\cite{piltner1985,qinhe2009}, for plate analysis~\cite{jirousekwroblewski1995,qin1995,choochoi2010}. Recently, the idea of employing local solutions over arbitrary finite elements have been investigated in~\cite{hofreitherlanger2010,weisser2011,rjasanowwei$ss$er2012,weisser2012}. However, its convergence properties and accuracy when applied to linear elasticity needs to be investigated.

\subsection{Approach} In this paper, hybrid-Trefftz arbitrary polygons will be formulated and its convergence properties and accuracy will be numerically studied with a few benchmark problems in the context of linear elasticity. An optimal number of T-complete functions are chosen based on the number of nodes of the polygon and degrees of freedom per node. The salient features of the approach are: (a) Only the boundary of the element is discretized with 1D finite elements and (b) Explicit form of the shape functions and special numerical integration scheme is not required to compute the stiffness matrix. 

\subsection{Outline} The paper commences with an overview of the governing equations for elasticity and the corresponding Galerkin form. Section \ref{htfemreview} introduces a hybrid Trefftz type approximation over arbitrary polygons. The efficiency, the accuracy and the convergence properties of the HTFEM is demonstrated with a few benchmark problems in Section \ref{numex}. The numerical results from the HTFEM are compared with the analytical results and with the polygonal FEM with Wachspress interpolants, followed by concluding remarks in the last section.

\section{Governing equations and weak form}
\label{fempoly}
\subsection{Governing equations and weak form}
For a 2D static elasticity problem defined in the domain $\Omega$ bounded by $\Gamma = \Gamma_u \bigcup \Gamma_t$, $\Gamma_u \bigcap \Gamma_t = \emptyset$, in the absence of body forces, the governing equation is given by:
\begin{equation}
\nabla_s \cdot \bvsig + \mathbf{b} = \mathbf{0} \hspace{0.25cm} \textup{in} \hspace{0.25cm}  \Omega
\label{eqn:geqn}
\end{equation}
with the following conditions prescribed on the boundary: 
\begin{align}
\uu &= \overline{\uu} \hspace{0.25cm} \textup{in} \hspace{0.25cm}  \Gamma_u \nonumber \\
\bvsig \cdot \mathbf{n} &= \overline{\mathbf{t}} \hspace{0.25cm} \textup{on} \hspace{0.25cm}  \Gamma_t
\end{align}
where $\bvsig$ is the stress tensor. The discrete equations for this problem are formulated from the Galerkin weak form:
\begin{equation}
\int\limits_{\Omega} \left( \nabla_s \uu \right)^{\rm T} \dd \left( \nabla_s \delta \uu \right)~\rmd \Omega - \int\limits_\Omega \left( \delta \uu \right)^{\rm T} \bb~\rmd \Omega - \int\limits_{\Gamma_t} \left( \delta \uu \right)^{\rm T} \overline{\mathbf{t}}~\rmd \Gamma = \mathbf{0}
\label{eqn:weakform}
\end{equation}
where $\uu$ and $\delta \uu$ are the trial and the test functions, respectively and $\dd$ is the material constitutive matrix. The FEM uses the following trial and test functions:
\begin{equation}
\uu^h(\xx) = \sum\limits_{I=1}^{NP} \bn_I(\xx) \mathbf{d}_I, \hspace{0.25cm} \delta\uu^h(\xx) = \sum\limits_{I=1}^{NP} \bn_I(\xx) \delta \mathbf{d}_I
\label{eqn:trialtestfn}
\end{equation}
where $NP$ is the total number of nodes in the mesh, $\bn$ is the shape function matrix and $\mathbf{d}_I$ is the vector of degrees of freedom associated with node $I$. Upon substituting \Eref{eqn:trialtestfn} into \Eref{eqn:weakform} and invoking the arbitrariness of  $\delta \uu$, we obtain the following discretized algebraic system of equations:
\begin{equation}
\kk \mathbf{d} = \mathbf{f}
\end{equation}
with
\begin{eqnarray} 
\kk &=& \int\limits_{\Omega^h} \bigb^{\rm T} \dd \bigb~\rmd \Omega \nonumber \\
\mathbf{f} &=&  \int\limits_{\Omega^h} \bn^{\rm T} \bb~\rmd \Omega + \int\limits_{\Gamma_t} \bn^{\rm T} \overline{\mathbf{t}}~\rmd \Gamma
\label{eqn:stiffmat}
\end{eqnarray}
where $\kk$ is the stiffness matrix and $\Omega^h$ is the discretized domain formed by the union of elements $\Omega^e$. The stiffness matrix is computed over each element and assembled to the global matrix. The size of the stiffness matrix depends on the number of nodes in an element.

\subsection{Generalization to arbitrary polygons}
\begin{figure}[htpb]
\centering
\subfigure[]{\includegraphics[scale=0.41]{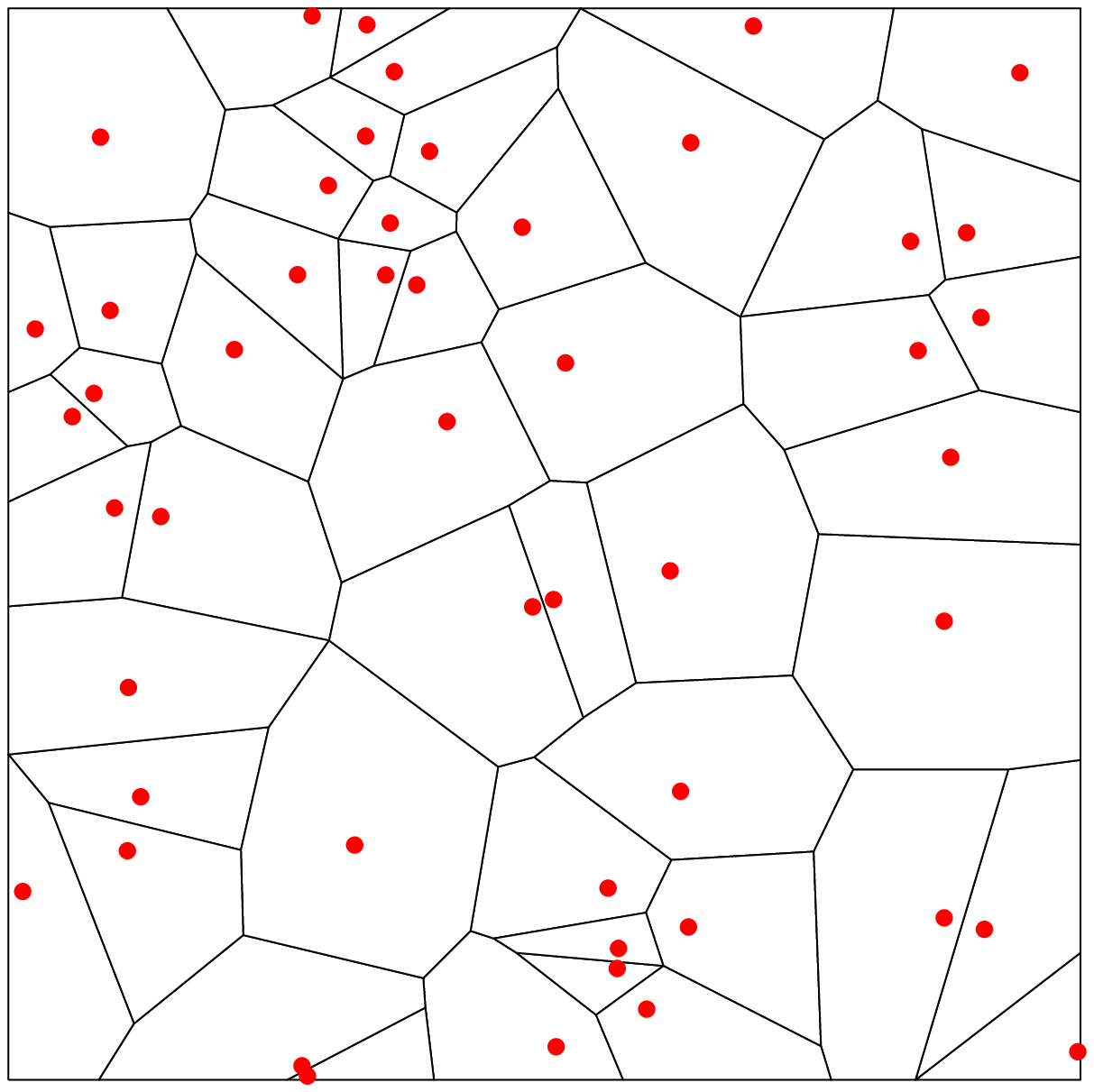}}
\subfigure[]{\includegraphics[scale=0.41]{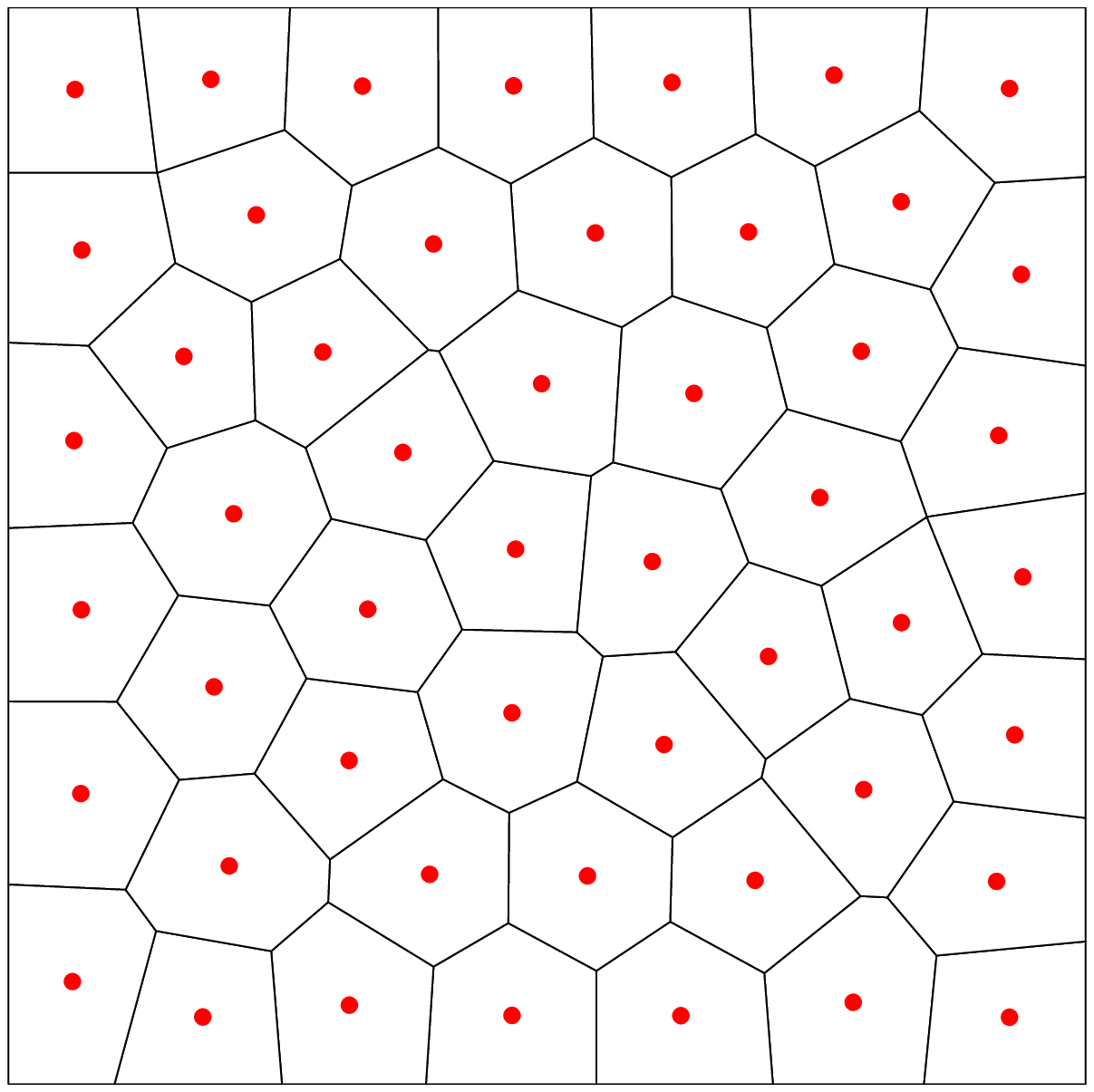}}
\caption{Vorono\"{i} tessellation of two set of data points, where the \textcolor{red}{red} dots are the seed points: (a) Scattered data set (b) Polygonal mesh after few iterations. }
\label{fig:vtessel}
\end{figure}

The growing interest in the generalization of FE over arbitrary meshes has opened up a new area of finite elements called `\emph{polygonal finite elements}'. In polygonal finite elements, the number of sides of an element is not restricted to three or four as in the case of 2D. The Vorono\"{i} tessellation is a fundamental geometrical construct to generate a polygonal mesh covering a given domain. Polygonal meshes can be generated from Vorono\"{i} diagrams. The Vorono\"{i} diagram is a subdivision of the domain into regions $V(p_I)$, such that any point in $V(p_I)$ is closer to node $p_I$ than to any other node. \fref{fig:vorfig} shows a Vorono\"{i} diagram of a point $P$. The first order Vorono\"{i} $V(N)$ is a subdivision of the Euclidean space $\mathbb{R}^2$ into convex regions, mathematically:
\begin{figure}
\centering
\includegraphics[scale=0.5]{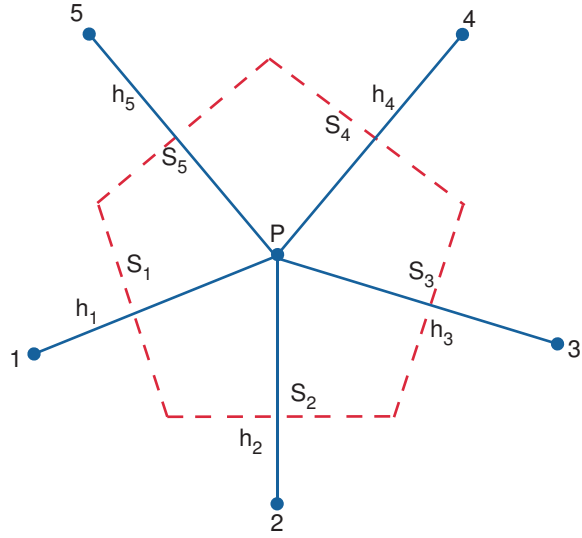}
\caption{Vorono\"{i} diagram of a point $P$.}
\label{fig:vorfig}
\end{figure}
\begin{equation}
T_I = \left\{ \xx \in \mathbb{R}^2 \colon d(\xx,\xx_I) < d(\xx,\xx_J) \forall J \neq I \right\}
\end{equation}
where $d(\xx_I,\xx_J)$, the Euclidean matrix, is the distance between $\xx_I$ and $\xx_J$. The quality of the generated polygonal mesh depends on the randomness in the scattered points. \fref{fig:vtessel} shows a typical Vorono\"{i} tessellation of two sets of scattered data set. The quality of a polygonal mesh determines the accuracy of the solution~\cite{sukumartabarraei2004}. To improve the quality of the Vorono\"{i} tessellation, the generating point of each Vorono\"{i} cell can be used as its center of mass, leading to a special type of Vorono\"{i} diagram, called the centroidal Vorno\"{i} tessellation (CVT)~\cite{duwang2005}. Sieger \textit{et al.,}~\cite{siegeralliez2010} presented an optimizing technique to improve the Vorono\"{i} diagrams for use in FE computations. 

\section{Overview of hybrid Trefftz finite element method}
\label{htfemreview} 
\begin{figure}[htpb]
\centering
\scalebox{0.8}{\input{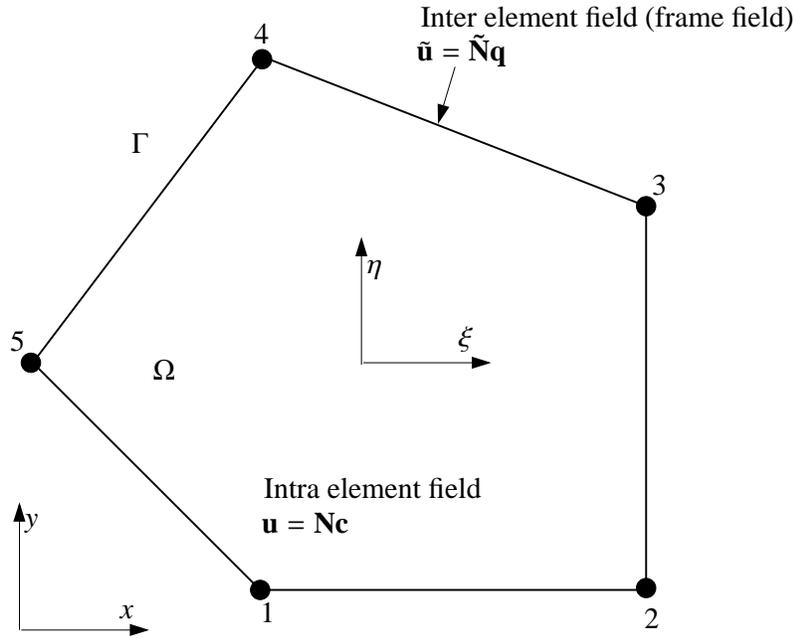}}
\caption{Hybrid Trefftz polygonal element with a description of inter and intra element field.}
\label{fig:htpfemidea}
\end{figure}

The basic idea in the hybrid Trefftz FEM is to employ the series of the homogeneous solution to the governing differential equation (see \Eref{eqn:geqn}) as the approximation function to model the displacement field within the domain and an independent set of functions to represent the boundary and to satisfy inter-element compatibility (see \fref{fig:htpfemidea}. The set of functions that are used to represent the dispalcement field within the domain are also called as T-complete set. The displacement field within an element can be written as:
\begin{equation}
\uu(\xx) = \sum\limits_{I} \bn_I(\xx) \mathbf{c}_I, \hspace{0.2cm} \xx \in \Omega
\label{eqn:intraele}
\end{equation}
where $\mathbf{c}$ are the vector of undetermined coefficients and $\mathbf{N}_I$ are the approximation functions that are selected from the series solution of the homogeneous part of the governing differential equation (see \Eref{eqn:geqn}). For linear elastostatics, based on the Mushelishvili's complex variable formulation, the $\mathbf{N}_I$ and the corresponding stress fields are given by~\cite{qin2005}:
\begin{eqnarray}
\bn_{Jk} &= \frac{1}{2G} \left\{ \begin{array}{c} \mathrm{Re}Z_{Jk} \\ \mathrm{Im} Z_{Jk}\end{array} \right\} \nonumber \\
\ttt_{Jk} &= \left\{ \begin{array}{c} \mathrm{Re}(R_{Jk} - S_{Jk}) \\ \mathrm{Re}(R_{Jk}+S_{Jk}) \\ \mathrm{Im}S_{Jk} \end{array} \right\}
\end{eqnarray}
where $J=1,2,3,4$, $k = 1,2,\cdots$, $Z_{1k} = i \kappa z^k + i \kappa z \overline{z}^{k-1};~~ Z_{2k} = \kappa z^k - k z \overline{z}^{k-1};~~Z_{3k} = i \overline{z}^k$ and $Z_{4k} = -\overline{z}^k$ and $R_{1k}=2ikz^{k-1}, R_{2k}=2kz^{k-1}, R_{3k} = 0, R_{4k} = 0$ and $S_{1k} = ik(k-1)z^{k-2}\overline{z}, S_{2k} = k(k-1)z^{k-2}\overline{z}, S_{3k} = iKz^{k-1}, S_{4k} = kz^{k-1}$.

\vspace{1pt}

However, the intra element displacement field given by \Eref{eqn:intraele} is non-conforming across the inter-element boundary. The unknown coefficients $\mathbf{c}$ are computed from the external boundary conditions and/or from the continuity conditions on the inter element boundary. Of the various methods available to enforce these conditions, in this study, we use the hybrid technique. In this technique, the elements are linked through an auxiliary conforming displacement frame which has the same form as in the conventional FEM. The displacement field on the element boundary, or otherwise called as frame is given by:
\begin{equation}
\tilde{\uu}(\xx) = \sum\limits_{I} \tilde{\bn}_I(\xx) \qq_I, \hspace{0.2cm} \xx \in \Gamma
\label{eqn:interele}
\end{equation}
where $\qq_I$ are the unknowns of the problem and $\tilde{\bn}_I$ are the standard 1D FE shape functions. To satisfy the inter-element continuity, a modified variational form is employed~\cite{qin2005}, given by:
\begin{equation}
\Pi = \frac{1}{2} \int\limits_\Omega \bvsig^{\rm T} \dd^{-1} \bvsig~\rmd \Omega - \int\limits_{\Gamma_{e1}} \mathbf{t} \delta \uu~\rmd \Gamma + \int\limits_{\Gamma_{e2}} (\overline{\mathbf{t}} - \mathbf{t}) \uu~\rmd \Gamma - \int\limits_{\Gamma_{eI}} \mathbf{t}\tilde{\delta \uu}~\rmd \Gamma
\label{eqn:modifyvariprin}
\end{equation}
where $\Gamma_{eI}$ is the inter-element boundary, $\Gamma_{e1} = \Gamma_u \cap \Gamma_e$ and $\Gamma_{e2} = \Gamma_\sigma \cap \Gamma_e$. The minimization of the modified variational principle given by \Eref{eqn:modifyvariprin} leads to the following system of algebraic equations:
\begin{equation}
\kk \mathbf{d} = \mathbf{f}
\end{equation}
where
\begin{eqnarray*}
\kk &= \mathbf{G}^{\rm T} \mathbf{H}^{-1} \mathbf{G} \nonumber \\
\ff  &= \int\limits_\Gamma \tilde{\bn} \overline{\mathbf{t}}~\rmd \Gamma
\end{eqnarray*}
and
\begin{equation*}
\mathbf{H} = \int\limits_\Gamma \mathbf{Q}^{\rm T} \bn~\rmd \Gamma; \hspace{0.25cm}
\mathbf{G} = \int\limits_\Gamma \mathbf{Q}^{\rm T} \tilde{\bn}~\rmd \Gamma .
\end{equation*}
where $\mathbf{Q} = \mathbf{A} \ttt$, and,
\begin{equation*}
\mathbf{A} = \left[ \begin{array}{ccc} n_1 & 0 & n_2 \\ 0 & n_2 & n_1 \end{array} \right] \hspace{0.2cm} \textup{and} \hspace{0.2cm} 
\ttt = \nabla \bn
\end{equation*}
where $n_1$ and $n_2$ are the outward normals. The integrals in the above equation can be computed by employing standard Gaussian quadrature rules. It is noted that in the computation of the stiffness matrix, we have to compute the inverse of the matrix $\mathbf{H}$. The necessary condition for the matrix $\mathbf{H}$ to be of full rank is
\begin{equation}
m_{\rm mim} = N_{\rm dof} - 1
\end{equation}
where $N_{\rm dof}$ is the total number of degrees of freedom of the element and $m_{\rm min}$ is the minimum number of T-complete functions to be used. Additionally, if $m_{\rm min}$ does not guarantee a matrix with full rank, full rank can be achieved by suitably increasing the number of T-complete functions.

\section{Numerical Examples}
\label{numex} In this section, we present the convergence and accuracy of the arbitrary polygons with local Trefftz functions using benchmark problems in the context of linear elasticity. The results from the proposed approach are compared with analytical solution where available and with the conventional polygonal finite element method with Laplace interpolants. To discuss the results, we employ the following convention:
\begin{itemize}
\item PFEM - Polygonal finite element method with Laplace/Wachspress interpolants (conventional approach). The numerical integration within each element is done by sub-dividing the polygon into triangles and employing a sixth order Dunavant quadrature rule.
\item HT-PFEM - hybrid Trefftz polygonal finite element method. Within each polygon, T-complete functions are employed to compute the stiffness matrix. One dimensional Gaussian quadrature is employed along the boundary of the polygon and the order of the quadrature depends on the number of T-complete functions employed.
\end{itemize}

The built-in Matlab\textsuperscript{\textregistered} function {\small voronoin} and Matlab\textsuperscript{\textregistered} functions in {\small PolyTop}~\cite{talischipaulino2012} for building the mesh-connectivity are used to create the polygonal meshes. For the purpose of error estimation, we employ the relative error in the displacement norm and in the energy norm, given by:
\paragraph{Displacement norm}
\begin{equation}
|| \uu - \uu^h||_{L^2(\Omega)} = \frac{ \sqrt{ \int\limits_{\Omega} (\uu-\uu^h) \cdot (\uu-\uu^h)~\rmd \Omega}}{\sqrt{ \int\limits_{\Omega} \uu \cdot \uu~\rmd \Omega}}
\end{equation}

\paragraph{Energy norm}
\begin{equation}
|| \uu - \uu^h||_{H^1(\Omega)} = \frac{ \sqrt{ \int\limits_{\Omega} (\bveps - \bveps^h)^{\rm T} \dd (\bveps - \bveps^h)~\rmd \Omega}}{\sqrt{ \int\limits_{\Omega} \bveps^{\rm T} \dd \bveps~\rmd \Omega}}
\end{equation}
where $\uu, \bveps$ are the analytical solution or a reference solution and $\uu^h, \bveps^h$ are the numerical solution.

\subsection{Cantilever beam bending}
\begin{figure}[htpb]
\centering
\scalebox{0.7}{\input{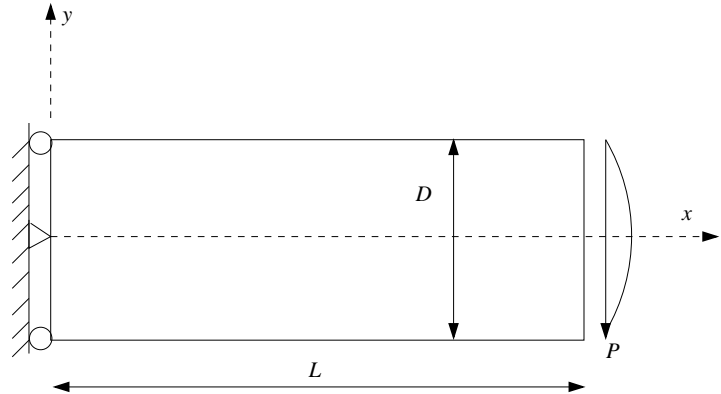}}
\caption{Cantilever beam: Geometry and boundary conditions.}
\label{fig:cantileverfig}
\end{figure}

\begin{figure}[htpb]
\centering
\includegraphics[scale=1]{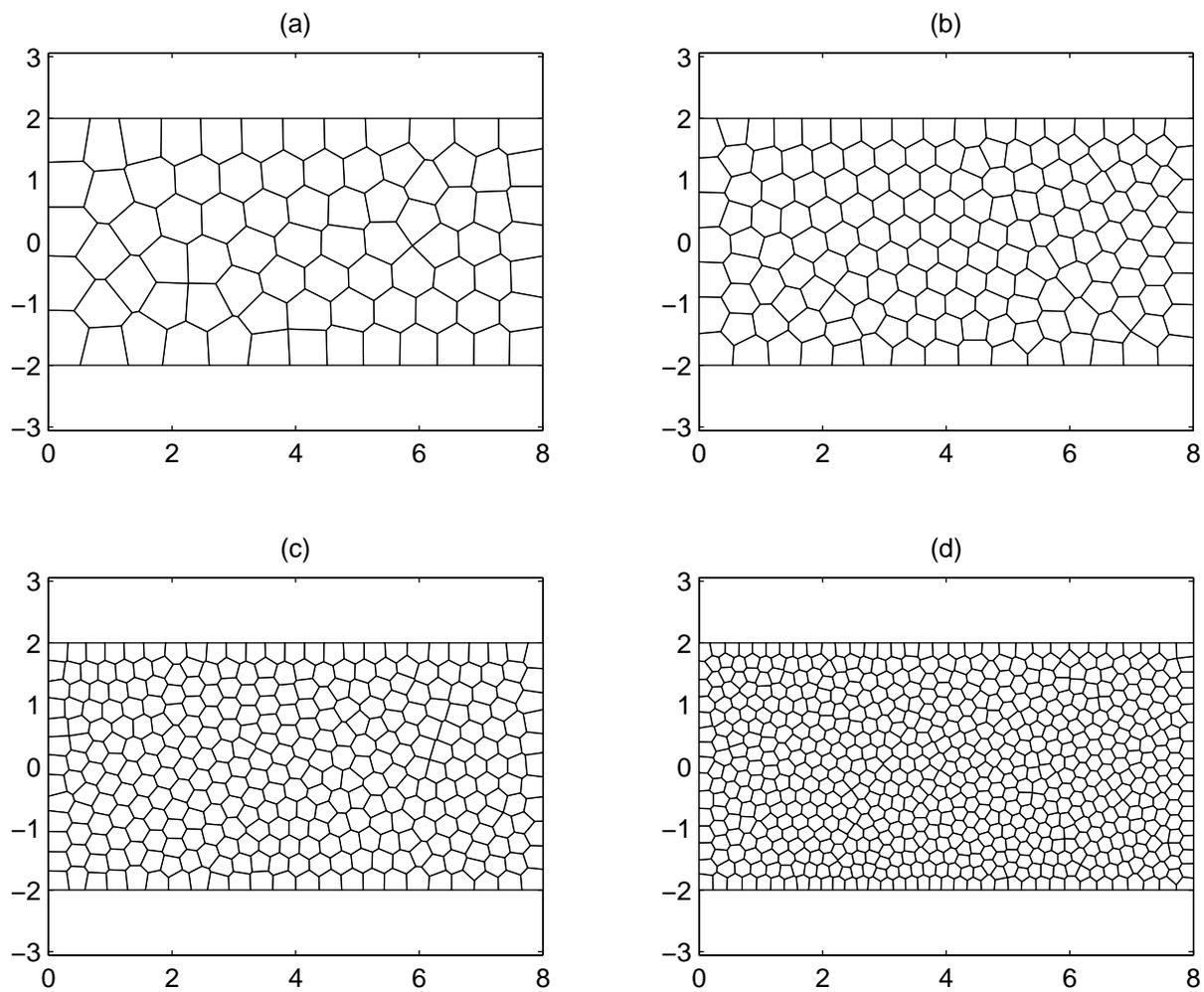}
\caption{Thick cantilever beam: domain discretized with polygonal elements: (a) 80 elements; (b) 160 elements; (c) 320 elements and (d) 640 elements.}
\label{fig:cantmesh}
\end{figure}

A two-dimensional cantilever beam subjected to a parabolic shear load at the free end is examined as shown in \fref{fig:cantileverfig}. The geometry is: length $L=$ 10, height $D=$ 2. The material properties are: Young's modulus, $E=$ 3e$^7$, Poisson's ratio $\nu=$ 0.25 and the parabolic shear force $P=$ 150. The exact solution for displacements are given by:

\begin{align}
u(x,y) &= \frac{P y}{6 \overline{E}I} \left[ (9L-3x)x + (2+\overline{\nu}) \left( y^2 - \frac{D^2}{4} \right) \right] \nonumber \\
v(x,y) &= -\frac{P}{6 \overline{E}I} \left[ 3\overline{\nu}y^2(L-x) + (4+5\overline{\nu}) \frac{D^2x}{4} + (3L-x)x^2 \right]
\label{eqn:cantisolution}
\end{align}
where $I = D^3/12$ is the moment of inertia, $\overline{E} = E$, $\overline{\nu} = \nu$ and $\overline{E} = E/(1-\nu^2)$, $\overline{\nu} = \nu/(1-\nu)$ for plane stress and plane strain, respectively. \fref{fig:cantmesh} shows a sample polygonal mesh used for this study.  

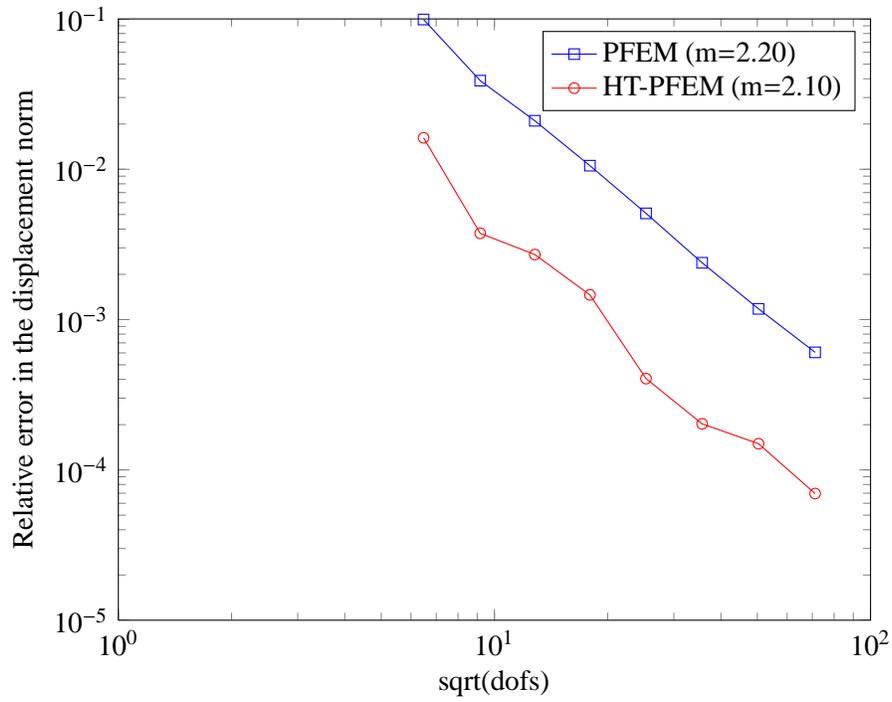
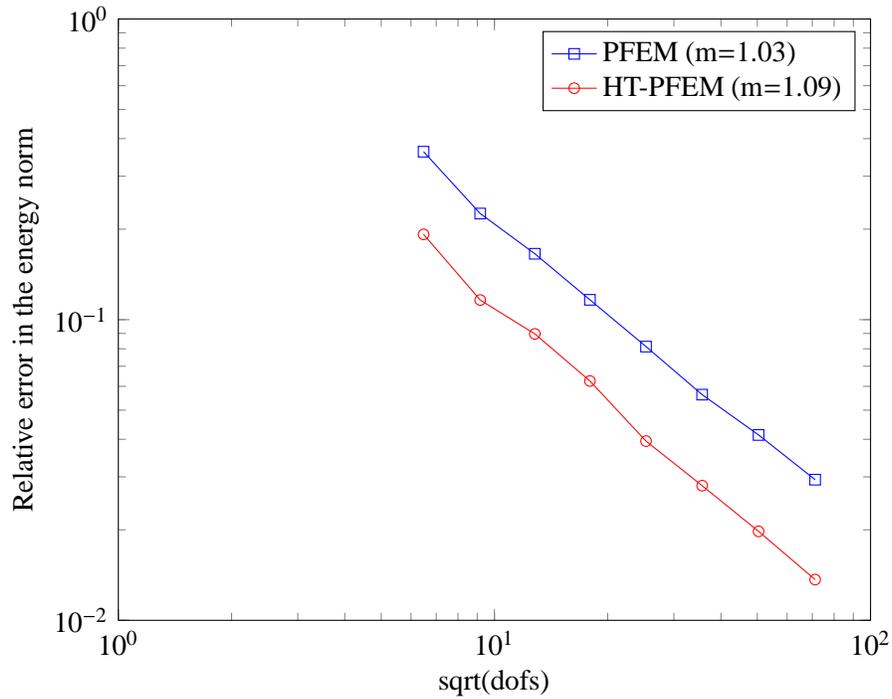
\begin{figure}
\centering
\newlength\figureheight 
\newlength\figurewidth 
\setlength\figureheight{8cm} 
\setlength\figurewidth{10cm}
\subfigure[]{
%
%
%
%
\begin{tikzpicture}

\begin{axis}[%
width=\figurewidth,
height=\figureheight,
scale only axis,
xmode=log,
xmin=1,
xmax=100,
xminorticks=true,
xlabel={sqrt(dofs)},
ymode=log,
ymin=1e-05,
ymax=0.1,
yminorticks=true,
ylabel={Relative error in the displacement norm},
legend style={draw=black,fill=white,legend cell align=left}
]
\addplot [
color=blue,
solid,
mark=square,
mark options={solid}
]
table[row sep=crcr]{
6.48074069840786 0.09913961\\
9.16515138991168 0.038849023\\
12.8062484748657 0.021016214\\
17.9443584449264 0.010565066\\
25.298221281347 0.005087769\\
35.665109000254 0.002386559\\
50.37856687124 0.001177957\\
71.2320152740325 0.0006053307587\\
};
\addlegendentry{PFEM (m=2.20)};

\addplot [
color=red,
solid,
mark=o,
mark options={solid}
]
table[row sep=crcr]{
6.48074069840786 0.016189419\\
9.16515138991168 0.003746207\\
12.8062484748657 0.002708105\\
17.9443584449264 0.001459322\\
25.298221281347 0.0004049331981\\
35.665109000254 0.0002024967679\\
50.37856687124 0.0001493225916\\
71.2320152740325 6.958269688e-05\\
};
\addlegendentry{HT-PFEM (m=2.10)};

\end{axis}
\end{tikzpicture}%
}
\subfigure[]{
%
%
%
%
\begin{tikzpicture}

\begin{axis}[%
width=\figurewidth,
height=\figureheight,
scale only axis,
xmode=log,
xmin=1,
xmax=100,
xminorticks=true,
xlabel={sqrt(dofs)},
ymode=log,
ymin=0.01,
ymax=1,
yminorticks=true,
ylabel={Relative error in the energy norm},
legend style={draw=black,fill=white,legend cell align=left}
]
\addplot [
color=blue,
solid,
mark=square,
mark options={solid}
]
table[row sep=crcr]{
6.48074069840786 0.362003885338265\\
9.16515138991168 0.225405503482058\\
12.8062484748657 0.165654803733547\\
17.9443584449264 0.116407169023218\\
25.298221281347 0.0813583861688517\\
35.665109000254 0.056314855944058\\
50.37856687124 0.0413500906891388\\
71.2320152740325 0.029359830520628\\
};
\addlegendentry{PFEM (m=1.03)};

\addplot [
color=red,
solid,
mark=o,
mark options={solid}
]
table[row sep=crcr]{
6.48074069840786 0.192202518193701\\
9.16515138991168 0.11622759138862\\
12.8062484748657 0.0896668444855734\\
17.9443584449264 0.062517885440888\\
25.298221281347 0.0394507794599803\\
35.665109000254 0.0280472043526623\\
50.37856687124 0.0197641942031543\\
71.2320152740325 0.0136586931146431\\
};
\addlegendentry{HT-PFEM (m=1.09)};

\end{axis}
\end{tikzpicture}%
}
\caption{Bending of thick cantilever beam: Convergence results for (a) the relative error in the displacement norm $(L^2)$ and (b) the relative error in the energy norm. The rate of convergence is also shown, where $m$ is the average slope. In case of the polygonal SBFEM, $p$ denotes the order of the shape functions along each edge of the polygon.}
\label{fig:cantiConveResults}
\end{figure}

The numerical convergence of the relative error in the displacement norm and the relative error in the energy norm is shown in \fref{fig:cantiConveResults}. The results from different approaches are compared with the available analytical solution. Both the Polygonal FEM and the HT-PFEM yields optimal convergence in $L^2$ and $H^1$ norm. It is seen that with mesh refinement, both the methods converge to the exact solution. An estimation of the convergence rate is also shown. From \fref{fig:cantiConveResults}, it can be observed that the HT-PFEM yields more accurate results. 

\subsection{Infinite plate with a circular hole}
\begin{figure}[htpb]
\centering
\scalebox{0.9}{\input{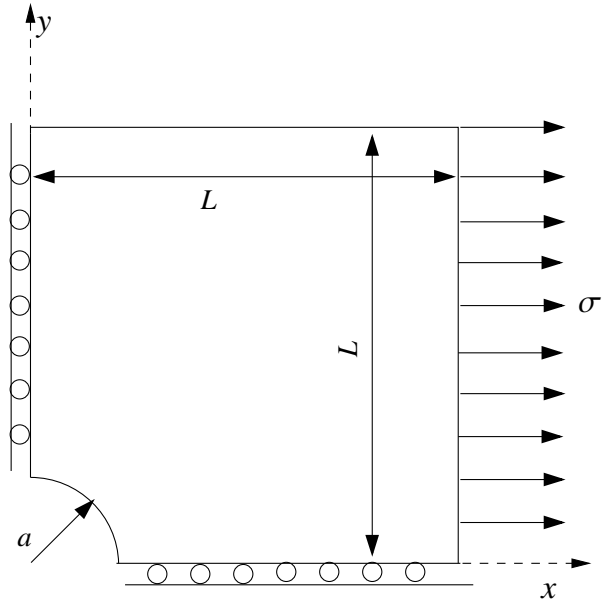}}
\caption{Infinite plate with a circular hole.}
\label{fig:twocrkhole}
\end{figure}

\begin{figure}[htpb]
\centering
\includegraphics[scale=1]{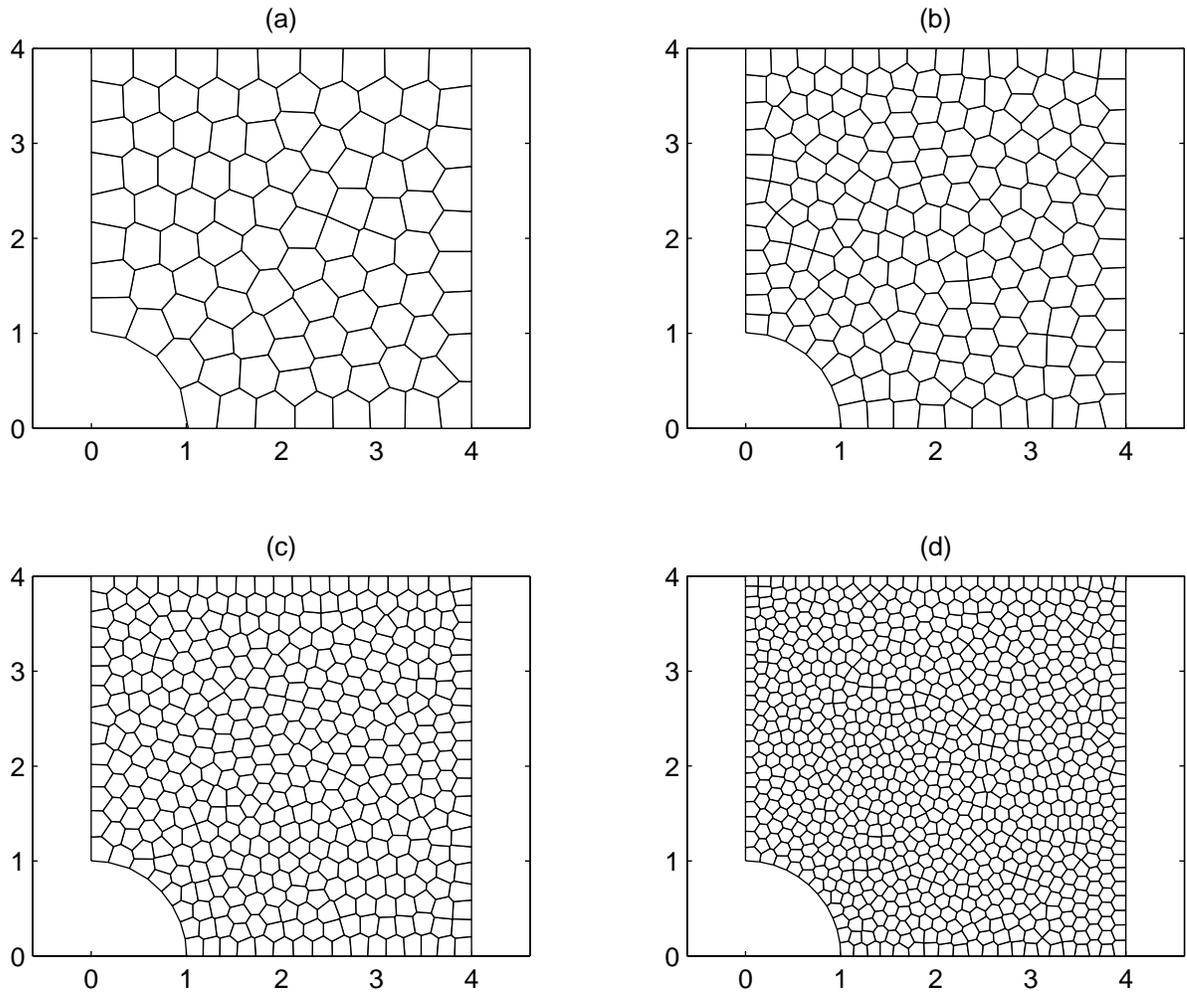}
\caption{Plate with a circular hole: domain discretized with polygonal elements: (a) 100 elements; (b) 200 elements; (c) 400 elements and (d) 800 elements.}
\label{fig:phwmesh}
\end{figure}

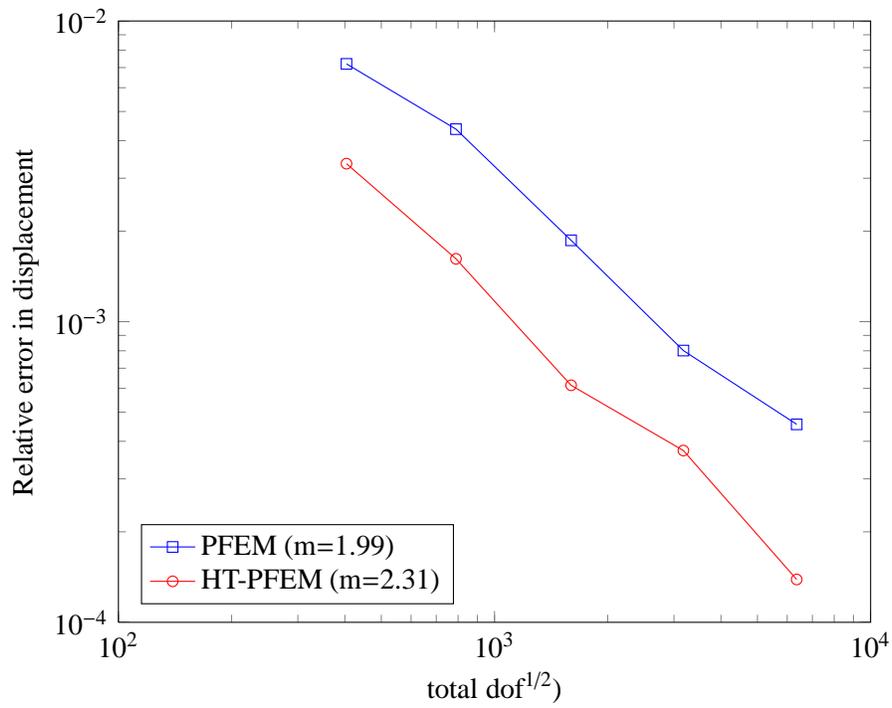
\begin{figure}
\centering
\setlength\figureheight{8cm} 
\setlength\figurewidth{10cm}
%
%
%
%
\begin{tikzpicture}

\begin{axis}[%
width=\figurewidth,
height=\figureheight,
scale only axis,
xmode=log,
xmin=100,
xmax=10000,
xminorticks=true,
xlabel={$\text{total dof}^{\text{1/2}}\text{)}$},
ymode=log,
ymin=0.0001,
ymax=0.01,
yminorticks=true,
ylabel={Relative error in displacement},
legend style={draw=black,fill=white,legend cell align=left},
legend pos = south west
]
\addplot [
color=blue,
solid,
mark=square,
mark options={solid}
]
table[row sep=crcr]{
404 0.007201360693307\\
790 0.004368235117637\\
1598 0.00186467533149\\
3178 0.0008014627252884\\
6360 0.000455433128974\\
};
\addlegendentry{PFEM (m=1.99)};

\addplot [
color=red,
solid,
mark=o,
mark options={solid}
]
table[row sep=crcr]{
404 0.003353954282388\\
790 0.001617946947749\\
1598 0.0006142561436902\\
3178 0.0003727080412345\\
6360 0.0001388314545548\\
};
\addlegendentry{HT-PFEM (m=2.31)};

\end{axis}
\end{tikzpicture}%
\caption{Infinite plate with a circular hole: Convergence results for the relative error in the displacement norm $(L^2)$. The rate of convergence is also shown, where $m$ is the average slope.}
\label{fig:plateHoleConveResults}
\end{figure}

In this example, consider an infinite plate with a traction free hole under uniaxial tension $(\sigma=$1$)$ along $x-$axis~\fref{fig:twocrkhole}. The exact solution of the principal stresses in polar coordinates $(r,\theta)$ is given by:
\begin{eqnarray}
\sigma_{11}(r,\theta) &= 1 - \frac{a^2}{r^2} \left( \frac{3}{2} (\cos 2\theta + \cos 4\theta) \right) + \frac{3a^4}{2r^4} \cos 4\theta \nonumber \\
\sigma_{22}(r,\theta) &= -\frac{a^2}{r^2} \left( \frac{1}{2}(\cos 2\theta - \cos 4\theta) \right) - \frac{3a^4}{2r^4} \cos 4\theta \nonumber \\
\sigma_{12}(r,\theta) &= -\frac{a^2}{r^2} \left( \frac{1}{2}(\sin 2\theta + \sin 4\theta) \right) + \frac{3a^4}{2r^4} \sin 4\theta
\end{eqnarray}
where $a$ is the radius of the hole. Owing to symmetry, only one quarter of the plate is modeled. \fref{fig:phwmesh} shows a typical polygonal mesh used for the study. The material properties are: Young's modulus $E=$ 10$^5$ and Poisson's ratio $\nu=$ 0.3. In this example, analytical tractions are applied on the boundary. The domain is discretized with polygonal elements and along each edge of the polygon, the shape function is linear. The convergence rate in terms of the displacement norm is shown in \fref{fig:plateHoleConveResults}. The relative error in the displacement norm for the PFEM and HT-PFEM is shown in \fref{fig:plateHoleConveResults}. It can be seen that the HT-PFEM yields more accurate results when compared with the PFEM. The HT-PFEM yields slightly a better convergence rate when compared to the PFEM. All the techniques converge to exact energy with mesh refinement. 

\subsection{Circular beam}
\begin{figure}[htpb]
\centering
\scalebox{0.7}{\input{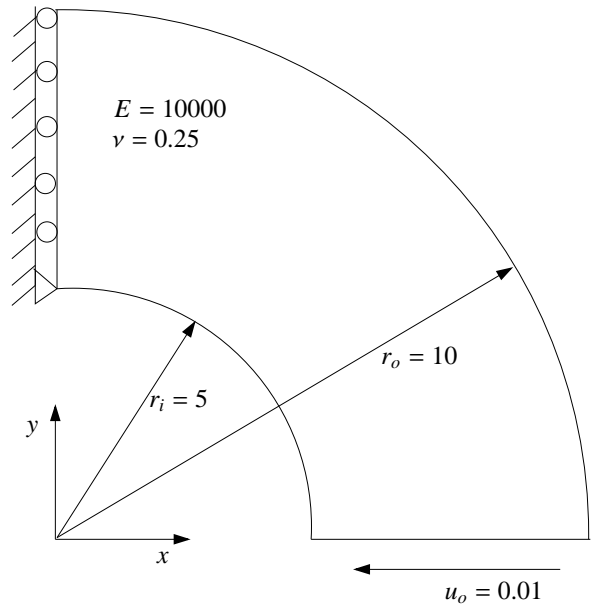}}
\caption{Curved cantilever circular beam: geometry, material property and boundary conditions.}
\label{fig:circbeam}
\end{figure}

\begin{figure}[htpb]
\centering
\includegraphics[scale=0.8]{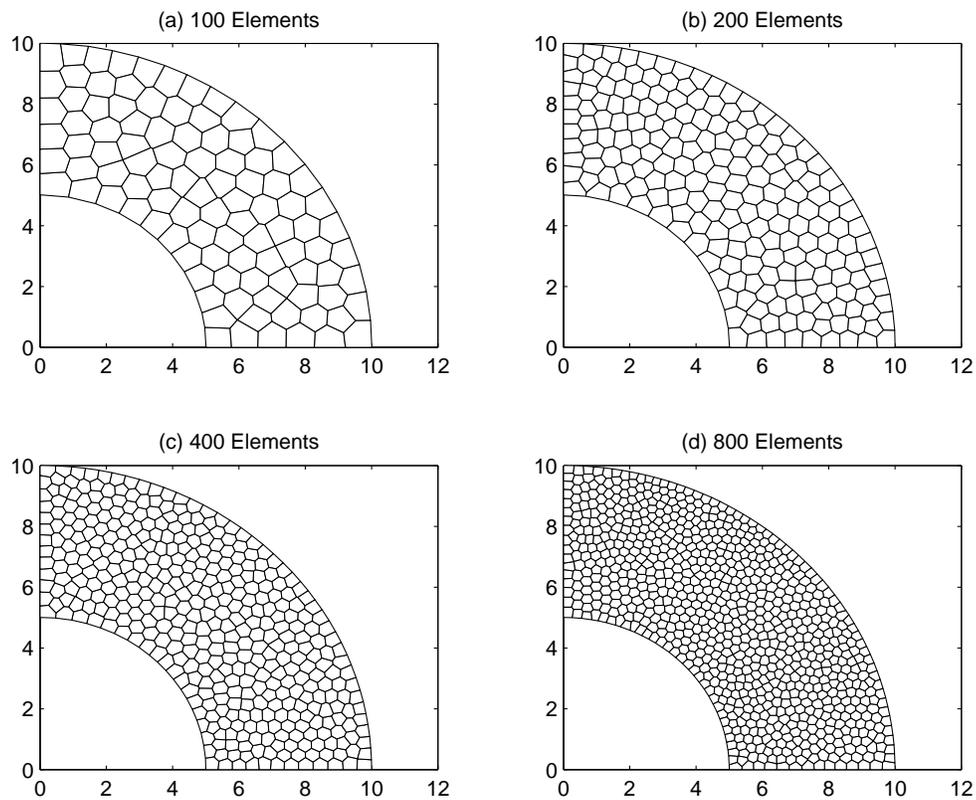}
\caption{Circular beam discretized with polygonal elements.}
\label{fig:cirmesh}
\end{figure}

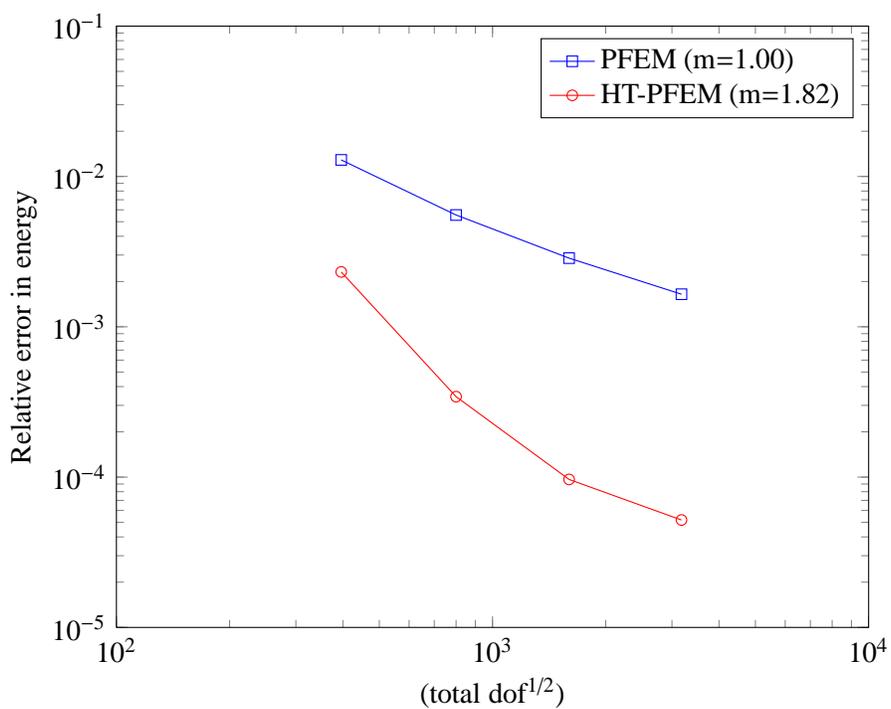
\begin{figure}
\centering
\setlength\figureheight{8cm} 
\setlength\figurewidth{10cm}
%
%
%
%
\begin{tikzpicture}

\begin{axis}[%
width=\figurewidth,
height=\figureheight,
scale only axis,
xmode=log,
xmin=100,
xmax=10000,
xminorticks=true,
xlabel={$\text{(total dof}^{\text{1/2}}\text{)}$},
ymode=log,
ymin=1e-05,
ymax=0.1,
yminorticks=true,
ylabel={Relative error in energy},
legend style={draw=black,fill=white,legend cell align=left}
]
\addplot [
color=blue,
solid,
mark=square,
mark options={solid}
]
table[row sep=crcr]{
396 0.0128842606958466\\
800 0.00554031147775727\\
1598 0.00286119115962027\\
3180 0.00164627075378143\\
};
\addlegendentry{PFEM (m=1.00)};

\addplot [
color=red,
solid,
mark=o,
mark options={solid}
]
table[row sep=crcr]{
396 0.0023170323270061\\
800 0.000342885436569265\\
1598 9.64908203462278e-05\\
3180 5.18450166141907e-05\\
};
\addlegendentry{HT-PFEM (m=1.82)};

\end{axis}
\end{tikzpicture}%
\caption{Circular beam: Convergence results for the relative error in the energy norm $(L^2)$. The rate of convergence is also shown, where $m$ is the average slope.}
\label{fig:cBeamConveResults}
\end{figure}

As a last example, consider a circular cantilevered beam as shown in \fref{fig:circbeam}. The beam is subjected to a prescribed displacement $u_o=$ -0.01 at the free end. The material property and boundary conditions considered for this study are shown in \fref{fig:circbeam}. The material is assumed to linear elastic and in a state of plane stress. The exact solution for the elastic energy is given by:
\begin{equation}
U = \frac{1}{\pi} ( \rm{ln}2 - 0.6)
\end{equation}
The circular beam is discretized with arbitrary polygons as shown in \fref{fig:cirmesh}. The convergence in the relative error in the energy norm is shown in \fref{fig:cBeamConveResults}. It can be seen that the HT-PFEM yields more accurate results when compared with the PFEM. The HT-PFEM yields a convergence rate of 1.82 and the PFEM yields a convergence rate of 1.0. Both the methods converge to the exact energy with mesh refinement.

\section{Concluding Remarks} In this paper, we constructed a hybrid T-Trefftz polygonal finite element by employing the T-complete set and an independent auxiliary field. The convergence and the accuracy of the proposed approach is demonstrated with a few benchmark examples. From the numerical studies presented, it is seen that the proposed approach yields more accurate results when compared to polygonal finite elements with Laplace interpolants. With hybrid Trefftz approach, special finite elements with embedded cracks/voids can be constructed. In the spirit of hybrid analytical XFEM~\cite{rethoreroux2010}, extended SBFEM~\cite{natarajansong2013} and hybrid crack element~\cite{xiaokarihaloo2007}, the region around the crack tip can be replaced with the special Trefftz element. The main advantage of such an approach would be the accurate representation of the crack tip singularities and the computation of the stress intensity factors directly. Another direction for future work is the study of hybrid Trefftz approach over arbitrary polyhedras.

\section*{Acknowledgements} Sundararajan Natarajan would like to acknowledge the financial support of the School of Civil and Environmental Engineering, The University of New South Wales for his research fellowship since September 2012.

\bibliographystyle{elsarticle-num}
\bibliography{htfem}

\end{document}